\documentclass[12pt,oneside,reqno]{amsart}
\usepackage{mathrsfs}
\usepackage{graphics}
\usepackage{amssymb}
\newcommand{\dif}{\mathrm{d}}

\newcommand{\be}{\begin{eqnarray}}
\newcommand{\ee}{\end{eqnarray}}
\newcommand{\ce}{\begin{eqnarray*}}
\newcommand{\de}{\end{eqnarray*}}
\newtheorem{theorem}{Theorem}[section]
\newtheorem{lemma}[theorem]{Lemma}
\newtheorem{remark}[theorem]{Remark}
\newtheorem{definition}[theorem]{Definition}
\newtheorem{proposition}[theorem]{Proposition}
\newtheorem{Example}[theorem]{Example}
\newtheorem{corollary}[theorem]{Corollary}

\def\[{{\Big[}}
\def\]{{\Big]}}
\def\<{{\langle}}
\def\>{{\rangle}}
\def\({{\Big(}}
\def\){{\Big)}}

\def\no{\nonumber}
\def\bt{\begin{theorem}}
\def\et{\end{theorem}}
\def\bl{\begin{lemma}}
\def\el{\end{lemma}}
\def\br{\begin{remark}}
\def\er{\end{remark}}
\def\bx{\begin{Example}}
\def\ex{\end{Example}}
\def\bd{\begin{definition}}
\def\ed{\end{definition}}
\def\bp{\begin{proposition}}
\def\ep{\end{proposition}}
\def\bc{\begin{corollary}}
\def\ec{\end{corollary}}

\def\cB{{\mathcal B}}
\def\cC{{\mathcal C}}

\def\cF{{\mathcal F}}

\def\cU{{\mathcal U}}

\def\mE{{\mathbb E}}

\def\mN{{\mathbb N}}

\def\mP{{\mathbb P}}

\def\mR{{\mathbb R}}

\def\mU{{\mathbb U}}

\def\mX{{\mathbb X}}

\def\geq{\geqslant}
\def\leq{\leqslant}

\begin{document}

\allowdisplaybreaks

\title{Topological equivalence for discontinuous random dynamical systems and applications*}

\author{Huijie Qiao$^1$ and Jinqiao Duan$^2$}



\thanks{*This work was supported by NSF of China (No. 11001051,  10971225 and 11028102), and by the NSF Grant  1025422.}


\date{\today}

\dedicatory{1. Department of Mathematics,
Southeast University\\
Nanjing, Jiangsu 211189,  China\\
hjqiaogean@yahoo.com.cn\\
2. Institute for Pure and Applied Mathematics (IPAM)\\University of California,
Los Angeles, CA 90095, USA\\ \& \\Department of Applied Mathematics, Illinois Institute of Technology\\
Chicago, IL 60616, USA\\
duan@iit.edu}

\begin{abstract}
After defining non-Gaussian  L\'evy processes for two-sided time, stochastic
differential equations with such L\'evy processes are considered. Solution paths for these stochastic differential equations have countable jump discontinuities in time. Topological equivalence (or conjugacy) for such an It\^o stochastic
differential equation and its transformed random differential equation  is established. Consequently, a stochastic Hartman-Grobman theorem is proved for the linearization of the It\^o stochastic differential
equation. Furthermore, for Marcus stochastic differential equations,
this topological equivalence   is  used to
prove existence of global random attractors.
\end{abstract}

\maketitle \rm

\section{Introduction}

Stochastic dynamical systems arise as mathematical models for
complex phenomena  under random fluctuations. They have been actively    studied when the fluctuations are Gaussian \cite{fs, il2, is}.
Non-Gaussian random fluctuations are, however, more widely observed in various
areas such as geophysics, biology, seismology, electrical engineering
and finance \cite{Woy, lb}. L\'evy processes are a class of non-Gaussian
processes whose sample paths are discontinuous in time.
For a dynamical system  driven by a L\'evy process, almost all paths or
orbits  have countable jump discontinuities in time.

  Discontinuous random dynamical systems or  cocycles (Definition \ref{mrds}) generated by stochastic
differential equations (SDEs) with L\'evy processes  have attracted attention more recently \cite{kuni, ks, q1, q2}.


In this paper, we consider topological equivalence between discontinuous cocycles,   generated either by SDEs with L\'evy processes or by related differential equations with random coefficients (i.e., random differential equations, or RDEs). Let us recall the  definition on topological equivalence or conjugacy \cite{la, il2}.

\bd\label{conj} Two random dynamical systems $\varphi$ and $\psi$ on
$\mR^d$ are called conjugate or topologically equivalent, if there exists a random homeomorphism
$H: \Omega\times\mR^d\mapsto\mR^d$ such that for all $(\omega,t)$,
$$
\psi_t(\omega,\cdot)=H(\theta_t\omega,\cdot)\circ\varphi_t(\omega,\cdot)\circ
H(\omega,\cdot)^{-1},
$$
where $H(\omega,\cdot)^{-1}$ stands for the inverse mapping of
$x\mapsto H(\omega,x)$. The homeomorphism $H$ is called a cohomology of
$\varphi$ and $\psi$. \ed

A cohomology of two random dynamical systems   is a random coordinate transformation, which transforms
the dynamical behavior for one of them into the dynamical behavior
for the other. Moreover, it does not change the intrinsic asymptotic
notions of these random dynamical systems, such as Lyapunov
exponents and random attractors \cite{il2}.

Conjugacy has been applied to study  SDEs with Brownian motion  in \cite{il2} and stochastic
partial differential equations with Brownian motion  in \cite{fl}.
It has also been used to examine the stochastic flows for SDEs with L\'evy processes in \cite{q1}. Consider the following
SDE in $\mR^d$: (Section \ref{conitosde})
\be\left\{\begin{array}{l}
\dif X_{t}=a\big(X_{t}\big)\,\dif t+\sigma_i\big(X_t\big)\dif L^i_{t}, \quad t\geq 0,\\
X_0=x,
\end{array}
\right. \label{Eq1}
\ee
and the RDE
\be\left\{\begin{array}{l}
\dif Y_t=\left(\frac{\partial\tilde{H}_t}{\partial x}\right)^{-1}(\omega,Y_t)a(\tilde{H}_t(\omega,Y_t))\dif t, \quad t\geq 0,\\
Y_0=x,
\end{array}
\right. \label{Eq2}
\ee
where the summation convention $\sum_i a_i b_i = a_i b_i$ is used, $L_t^i$'s are L\'evy processes, and  $\tilde{H}_t(\omega,x)$ is the
solution of the equation
\be
\tilde{H}_t(x)=x+\int^t_{0}\sigma_i(\tilde{H}_s(x))\dif L^i_{s}.
\label{Eq3}
\ee
In \cite{q1}, under appropriate conditions, the
first named author transformed SDE (\ref{Eq1}) to
 RDE  (\ref{Eq2}) by $\tilde{H}_t(\omega,x)$.
  Then a homeomorphism stochastic flow property for Eq.(\ref{Eq2}) is proved in order to get the
homeomorphism stochastic flow property for Eq.(\ref{Eq1}). Since
$\tilde{H}_t(\omega,x)$ is generally different from
$\tilde{H}_0(\theta_t\omega,x)$ for $t>0$ and $\omega\in\Omega$, we need to modify SDE (\ref{Eq3}) to
find a cohomology for topological equivalence.

In this paper, we transform   SDEs with L\'evy processes to RDEs by a cohomology. Because
pathwise arguments  for solutions of RDEs are more readily available than that of SDEs,
dynamical problems such as linearization (i.e., Hartman-Grobman theorem) and random attractors of SDEs can be dealt with.
Due to discontinuity in time for the solution paths of the   SDEs with L\'evy processes, the differentiation operation is conducted
via Fubini theorem and integration by parts, instead of the It\^o-ventzell formula. Additionally, for Marcus SDEs, we could only
construct these cohomologies in special cases.

It is worth mentioning that a major source of discontinuous cocycles
are solution mappings, after a perfection procedure \cite{ks}, of
stochastic differential equations with L\'evy processes. To consider
these cocycles with two-sided time $\mR$, we need to define L\'evy
processes for two-sided time as well.


\medskip

This paper is arranged as follows. In Section \ref{prelim}, we
introduce discontinuous cocycles, and L\'evy processes for two-sided time. In Section \ref{itosde}, we study
conjugacy for It\^o SDEs and RDEs and apply the result to prove a stochastic Hartman-Grobman theorem
  for  SDEs with L\'evy processes. Conjugacy for Marcus SDEs and RDEs, and
existence of global attractors are studied in Section
\ref{marcussde}.   For readers' convenience, the It\^o-Ventzell
formula for c\`adl\`ag processes is placed in the Appendix, Section \ref{itit}.

The following conventions will be used throughout the paper: $C$
with or without indices will denote different positive constants
(depending on the indices) whose values may change from one place to
another.

\section{Preliminaries}\label{prelim}

In this section, we recall several basic concepts and results which will be
needed throughout the paper.

\subsection{Basic notations}
The usual scalar product and norm (or length) in $\mR^d$ are  denoted by $\<\cdot,\cdot\>$ and $|\cdot|$, respectively. Moreover, $\nabla$
is   the gradient of vector fields on $\mR^d$, and  $[\cdot,\cdot]$
is the Lie bracket  operation on vector fields.

Denote by $\cC_b^{m,\gamma}$ the set of functions
$f:\mR^d\mapsto\mR^d$ satisfying \ce
\sup\limits_{x\in\mR^d}\frac{|f(x)|}{1+|x|}+\sum\limits_{|\beta|=1}^m\sup\limits_{x\in\mR^d}
|D^\beta f(x)|+\sum\limits_{|\beta|=m}\sup\limits_{x,y\in\mR^d,x\neq
y} \frac{|D^\beta f(x)-D^\beta f(y)|}{|x-y|^\gamma}<\infty, \de
where $D^\beta f(x)$ stands for the $\beta$-order partial derivative
of $f(x)$ for $|\beta|\leq m, m\in\mN$ and $0<\gamma<1$.

\subsection{Probability space}

Let $D^*(\mR,\mR^d)$ be the set of all functions which are
c\`adl\`ag (right continuous with left limit at each time) for $t\geq 0$ and c\`agl\`ad (left continuous with right limit at each time) for $t\leq 0$, and take
values in $\mR^d$. We take canonical sample space $\Omega \triangleq  D^*(\mR,\mR^d)$. It can be made a
complete and separable metric space  with endowed  Skorohod
metric $\rho$ as follows (\cite{Billingsley, hwy}):
\ce
\rho(x,y):=\inf\limits_{\lambda\in\Lambda}\left\{\sup\limits_{s\neq
t}\left|\log\frac{\lambda(t)-\lambda(s)}{t-s}\right|
+\sum\limits_{m=1}^\infty\frac{1}{2^m}\min\Big\{1,
\rho_m^\circ(x^m,y^m)\Big\}\right\}
\de
for all $x,y\in\Omega$,
where $x^m(t):=g_m(t)x(t)$, $y^m(t):=g_m(t)y(t)$ with \ce
g_m(t):=\left\{\begin{array}{c}
1, ~\quad\mbox{if}~ |t|\leq m,\quad\\
m+1-|t|,\quad\mbox{if}~m<|t|<m+1,\quad\\
\quad 0, \qquad \mbox{if}~|t|\geq m+1,
\end{array}
\right. \de and \ce \rho_m^\circ(x,y):=\sup\limits_{|t|\leq
m}\left|x(t)-y(\lambda(t))\right|. \de Here $\Lambda$ denotes the
set of strictly increasing and continuous functions from $\mR$ to
$\mR$. We identify a function $\omega(t)$  with a (canonical) sample
path $\omega$ in the sample space $\Omega$.

The Borel $\sigma$-algebra  in the sample space $\Omega$  under the
topology induced by the Skorohod metric $\rho$ is denoted as $\cF$.
Note that $\cF=\sigma(\omega(t),t\in\mR)$. Let
$\mP$ be the unique probability measure which makes the canonical
process a L\'evy process for $t\in\mR$ (see Definition \ref{levy1}
and \ref{levy2} below). And we have the complete natural filtration
$\cF_s^t:=\sigma(\omega(u): s\leq u\leq t)\vee\mathcal{N}$ for
$s\leq t$ with respect to $\mP$.

\subsection{Discontinuous random dynamical systems}(\cite{la})

Define  the Wiener shift
\ce
(\theta_t\omega)(\cdot)=\omega(t+\cdot)-\omega(t), \quad
t\in\mR, \;\; \omega\in\Omega.
\de
Then $\{\theta_t\}$ is a one-parameter group on $\Omega$. In fact, $\Omega$ is invariant with respect to
$\{\theta_t\}$, i.e.
$$
\theta_t^{-1}\Omega=\Omega, \quad ~\mbox{for~all}~ t\in\mR,
$$
and $\mP$ is $\{\theta_t\}$-invariant, i.e.
$$
\mP(\theta_t^{-1}(B))=\mP(B), \qquad ~\mbox{for~ all}~ B\in\cF,
t\in\mR.
$$
Thus $(\Omega,\cF,\mP,(\theta_t)_{t\in\mR})$ is a metric dynamical
system (DS) and ergodic, i.e., all measurable
$\{\theta_t\}$-invariant sets have probability $0$ or $1$.

\bd\label{mrds} Let $(\mX,\cB)$ be a measurable space. A mapping
\ce
\varphi: \mR\times\Omega\times\mX\mapsto\mX, \quad
(t,\omega,x)\mapsto\varphi_t(\omega,x)
\de
with the following properties is called a measurable random dynamical system (RDS), or in short,  a discontinuous cocycle:

(i) Measurability: $\varphi$ is
$\cB(\mR)\otimes\cF\otimes\cB/\cB$-measurable,

(ii) Discontinuous cocycle (over $\theta$) property:
$\varphi(t,\omega)$ is c\`adl\`ag for $t\geq0$
and c\`agl\`ad for $t\leq0$  and furthermore
\be
\varphi_0(\omega,\cdot)&=&id_{\mX}, \quad\qquad ~\mbox{for~ all}~\omega\in\Omega,\label{perfect coc1}\\
\varphi_{t+s}(\omega,\cdot)&=&\varphi_t(\theta_s\omega,\cdot)\circ\varphi_s(\omega,\cdot),
~\mbox{for all}~ s,t\in\mR, \quad \omega \in \Omega.
\label{perfect coc2}
\ee
\ed

\subsection{Random attractors}\label{ranattr}

We recall the definition of a random attractor  (\cite{fl}) and a
theorem for its existence (\cite{fs}). Let $\varphi_t$ be a
discontinuous cocycle.

\bd\label{temset} A random bounded set $B(\omega)\subset\mR^d$ is
called tempered if \ce \lim\limits_{t\rightarrow+\infty}e^{-\beta
t}d(B(\theta_{-t}\omega))=0, ~\mbox{ for any }~ \beta>0, \de where
$d(A)=\sup\limits_{x\in A}|x|$. \ed

\bd\label{absorbset} A random set $K(\omega)$ is said to be an
absorbing set if for all tempered random bounded set
$B(\omega)\subset\mR^d$ there exists $t(\omega,B)>0$ such that \ce
\varphi_t(\theta_{-t}\omega,B(\theta_{-t}\omega))\subset K(\omega),
~\mbox{ for all }~ t\geq t(\omega,B). \de \ed

\bd\label{globattr} A random compact set $A(\omega)\subset\mR^d$ is
called a global random attractor of $\varphi$ if for all
$\omega\in\Omega$ \ce
\varphi_t(\omega,A(\omega))=A(\theta_t\omega), ~\mbox{ for all }~ t\geq 0,\\
\lim\limits_{t\rightarrow+\infty}dist(\varphi_t(\theta_{-t}\omega,B(\theta_{-t}\omega)),A(\omega))=0,
\de for all tempered random bounded set $B(\omega)\subset\mR^d$,
where $dist$ denotes the semi-Hausdorff distance \ce
dist(A,B)=\sup\limits_{x\in A}\inf\limits_{y\in B}|x-y|. \de \ed

Note that a random attractor is unique. The following existence
theorem is from \cite{fs}.

\bt\label{exiattr} If there exists a compact absorbing set
$K(\omega)$, then there exists a random attractor of
$\varphi$.
\et

\subsection{L\'evy processes}\label{levypro}

\subsubsection{L\'evy processes for $t\geq0$}\label{levypro1}

\bd\label{levy1} A stochastic process $L=(L_t)_{t\geq0}$ with $L_0=0$ a.s. is a
$d$-dimensional L\'evy process if

(i) $L$ has independent increments; that is, $L_t-L_s$ is
independent of $L_v-L_u$ if $(u,v)\cap(s,t)=\emptyset$;

(ii) $L$ has stationary increments; that is, $L_t-L_s$ has the same
distribution as $L_v-L_u$ if $t-s=v-u>0$;

(iii) $L_t$ is right continuous with left limit. \ed

Its characteristic function is given by \ce
\mE\big(\exp\{i\<z,L_t\>\}\big)=\exp\{t\Psi(z)\}, \quad z\in\mR^d.
\de The function $\Psi: \mR^d\rightarrow\mathcal {C}$ is called the
characteristic exponent of the L\'evy process $L$. By the
L\'evy-Khintchine formula, there exist a nonnegative-definite
$d\times d$ matrix $Q$, $b\in\mR^d$ and a measure $\nu$ on $\mR^d$
satisfying \be \nu(\{0\})=0 ~\mbox{and}~
\int_{\mR^d}(|u|^2\wedge1)\nu(\dif u)<\infty, \label{lemc} \ee such
that \be
\Psi(z)&=&-\frac{1}{2}\<z,Qz\>+i\<z,b\>\no\\
&&+\int_{\mR^d}\big(e^{i\<z,u\>}-1-i\<z,u\>1_{|u|\leq\delta}\big)\nu(\dif
u), \label{lkf} \ee where $\delta>0$ is a constant. $\nu$ is called
the L\'evy measure.

Set $\kappa_t:=L_t-L_{t-}$. Then $\kappa$ defines a stationary
$(\cF_0^t)_{t\geq0}$-adapted Poisson point process with values in
$\mR^d\setminus\{0\}$ and the characteristic measure $\nu$ (\cite{iw}). Let $N_{\kappa}((0,t],\dif u)$ be the counting measure
of $\kappa_{t}$, i.e., for $B\in\cB(\mR^d\setminus\{0\})$
$$
N_{\kappa}((0,t],B):=\#\{0<s\leq t: \kappa_s\in B\},
$$
where $\#$ denotes the cardinality of a set. The compensator measure
of $N_{\kappa}$ is given by
$$
\tilde{N}_{\kappa}((0,t],\dif u):=N_{\kappa}((0,t],\dif u)-t\nu(\dif
u).
$$
The L\'evy-It\^o theorem states that there exist a $d'$-dimensional
$(\cF_0^t)_{t\geq0}$-Brownian motion $W_t$, $0\leq d'\leq d$ and  a
$d\times d'$ matrix $A$ such that $L$ can be represented as \be
L_t&=&bt+AW_t+\int_0^t\int_{|u|\leq\delta}u\tilde{N}_{\kappa}(\dif s, \dif u)\no\\
&&+\int_0^t\int_{|u|>\delta}uN_{\kappa}(\dif s, \dif u).
\label{lif1} \ee

\subsubsection{L\'evy processes for $t\leq0$}\label{levypro2}

\bd\label{levy2} A stochastic process $L^-=(L^-_t)_{t\leq0}$ with $L^-_0=0$
a.s. is a $d$-dimensional L\'evy process if

(i) $L^-$ has independent increments; that is, $L^-_t-L^-_s$ is
independent of $L^-_v-L^-_u$ if $(v,u)\cap(t,s)=\emptyset$;

(ii) $L^-$ has stationary increments; that is, $L^-_t-L^-_s$ has the
same distribution as $L^-_v-L^-_u$ if $t-s=v-u<0$;

(iii) $L^-_t$ is left continuous with right limit. \ed

Its characteristic function, characteristic exponent and the
L\'evy-Khintchine formula are the same as those for the L\'evy
process $L=(L_t)_{t\geq0}$.

Set $\kappa^-_t:=L^-_t-L^-_{t+}$. Then $\kappa^-$ defines a
stationary $(\cF_t^0)_{t\leq0}$-adapted Poisson point process with
values in $\mR^d\setminus\{0\}$ and characteristic measure
$\nu^-$ (\cite{iw}). Define \ce N_{\kappa^-}([t,0),B):=\#\{t\leq
s<0: \kappa^-_s\in B\}, \de for $B\in\cB(\mR^d\setminus\{0\})$. The
compensator measure of $N_{\kappa^-}$ is given by
$$
\tilde{N}_{\kappa^-}([t,0),\dif u):=N_{\kappa^-}([t,0),\dif
u)+t\nu^-(\dif u),
$$
where $\nu^-$ is the L\'evy measure for $L^-$. By the similar proof
to that of the L\'evy-It\^o theorem in \cite{sa} we obtain that
there exist $b^-\in\mR^d$, a $d^{''}$-dimensional
$(\cF_t^0)_{t\leq0}$-Brownian motion $W^-_t$, $0\leq d^{''}\leq d$
and a $d\times d^{''}$ matrix $A^-$ such that $L^-$ can be
represented as \be
L^-_t&=&-b^-t-A^-W^-_t-\int_t^0\int_{|u|\leq\delta}u\tilde{N}_{\kappa^-}(\dif s, \dif u)\no\\
&&-\int_t^0\int_{|u|>\delta}uN_{\kappa^-}(\dif s, \dif u).
\label{lif2} \ee

\medskip

A L\'evy process with two-sided time, $t\in  \mathbb{R}$, is defined to satisfy both Definitions \ref{levy1} and \ref{levy2}.

\section{Conjugacy for It\^o SDE and RDE, and linearization of SDEs}\label{itosde}

In this section, we study conjugacy for It\^o SDEs and RDEs and
apply the result to study linearization of SDEs with L\'evy processes.

\subsection{Conjugacy for It\^o SDEs and RDEs}\label{conitosde}

Take a one-sided $m$-dimensional L\'evy process
$$
L_t=bt+AW_t+\int^t_0\int_{|u|\leq\delta}u\,\tilde{N}_{\kappa}(\dif
s,\dif u), \qquad t\geq0,
$$
where $W_t$ is a $m'$-dimensional $(\cF_0^t)_{t\geq0}$-Brownian
motion, $0\leq m'\leq m$, $A=(a_{ij})$ is a $m\times m'$ matrix,
$0<\delta<1$, and consider the following It\^o SDE on $\mR^d$ for
$y\in\mR^d, \eta\in\mR$: \be\left\{\begin{array}{l}
\dif Y_t=\eta\sigma_i(Y_t)\dif L^i_t, \qquad t\geq0,\\
Y_0=y,
\end{array}
\right. \label{diffeq}
\ee
where the differential is in the It\^o sense (\cite{iw}).

\bp\label{diffeom} Assume that $\sigma_i(x)\in\cC_b^{1,\gamma}$ for
$i=1,2,\cdots,m$. If \ce \left(\begin{array}{c} y\\\eta \end{array}
\right)\mapsto\left(\begin{array}{c} y\\\eta \end{array}
\right)+\left(\begin{array}{c} \eta\sigma_i(y)u^i\\0 \end{array}
\right) \de is homeomorphic for any $u\in\{u\in\mR^m:
|u|\leq\delta\}$ and the Jacobian matrix \ce
I+\left(\begin{array}{c} \eta\frac{\partial(\sigma_i(y))^1}{\partial
y_1}u^i
\quad \eta\frac{\partial(\sigma_i(y))^1}{\partial y_2}u^i ~\cdots ~ \eta\frac{\partial(\sigma_i(y))^1}{\partial y_d}u^i \quad (\sigma_i(y))^1u^i\\
\eta\frac{\partial(\sigma_i(y))^2}{\partial y_1}u^i
\quad \eta\frac{\partial(\sigma_i(y))^2}{\partial y_2}u^i ~\cdots ~ \eta\frac{\partial(\sigma_i(y))^2}{\partial y_d}u^i \quad (\sigma_i(y))^2u^i\\
\vdots ~\qquad\qquad \vdots \qquad \cdots \qquad \vdots \qquad\qquad \vdots\\
\eta\frac{\partial(\sigma_i(y))^d}{\partial y_1}u^i
\quad \eta\frac{\partial(\sigma_i(y))^d}{\partial y_2}u^i ~\cdots ~ \eta\frac{\partial(\sigma_i(y))^d}{\partial y_d}u^i \quad (\sigma_i(y))^du^i\\
0 ~\qquad\qquad 0 \qquad \cdots \qquad 0 \qquad\qquad 0
\end{array}
\right) \de is invertible for any $y,\eta$ a.e. and
$u\in\{u\in\mR^m: |u|\leq\delta\}$, then the solution $Y_t^{y,\eta}$
to (\ref{diffeq}) defines a stochastic flow of $\cC^1$-diffeomorphisms. Moreover, $Y_t^{y,\eta}$ is differentiable in
$\eta$. \ep

\begin{proof}
For any $\eta\in\mR$, define $\eta_t:=\eta$ and rewrite
(\ref{diffeq}) as \be
\bar{Y}_t(\bar{y})&=&\bar{y}+\int^t_0F(\bar{Y}_s(\bar{y}))\dif s+\int^t_0G(\bar{Y}_s(\bar{y}))\dif W_s\no\\
&&+\int^t_0\int_{|u|\leq\delta}J(\bar{Y}_s(\bar{y}),u)\tilde{N}_{\kappa}(\dif
s,\dif u), \label{enlae} \ee where \ce
&&\bar{Y}_t=\left(\begin{array}{c}
Y_t^{y,\eta}\\
\eta_t
\end{array}
\right), \quad \bar{y}=\left(\begin{array}{c} y\\\eta \end{array}
\right),\\
&&F(\bar{Y}_t)=\left(\begin{array}{c}\eta\sigma_i(Y_t^{y,\eta})b^i\\
0\end{array}\right),
\quad J(\bar{Y}_t,u)=\left(\begin{array}{c}\eta\sigma_i(Y_t^{y,\eta})u^i\\ 0\end{array}\right),\\
&&G(\bar{Y}_t)=\left(\begin{array}{c}\eta\sigma_i(Y_t^{y,\eta})a^{i1}
\quad \eta\sigma_i(Y_t^{y,\eta})a^{i2} \quad \cdots \quad
\eta\sigma_i(Y_t^{y,\eta})a^{im'}\\0\qquad\qquad\qquad 0
\qquad\cdots \qquad\quad 0\end{array}\right). \de By Theorem 3.11 in
\cite{kuni}, for almost all $\omega$ and $t\geq0$,
$\bar{y}\mapsto\bar{Y}_t(\bar{y})$ is a stochastic flow of local
$\cC^1$-diffeomorphisms on $\mR^d\times\mR$.

 Now we prove that these diffeomorphisms are global. For any $R>0$,
$|\bar{y}|\leq R$. Fix $T>0$ and $p>d+1$. To (\ref{enlae}), by BDG
inequality in \cite[Theorem 2.11]{kuni} and H\"older's
inequality we have that \ce \mE\left(\sup\limits_{0\leq t\leq
T}|\bar{Y}_t|^p\right)&\leq&4^{p-1}|\bar{y}|^p
+4^{p-1}T^{p-1}\mE\left(\int^T_0|F(\bar{Y}_s(\bar{y}))|^p\dif s\right)\\
&&+4^{p-1}C\mE\left(\int^T_0\int_{|u|\leq\delta}|J(\bar{Y}_s(\bar{y}),u)|^p\nu(\dif u)\dif s\right)\\
&&+4^{p-1}C\mE\left(\int^T_0\int_{|u|\leq\delta}|J(\bar{Y}_s(\bar{y}),u)|^2\nu(\dif u)\dif s\right)^{\frac{p}{2}}\\
&&+4^{p-1}C\mE\left(\int^T_0|G(\bar{Y}_s(\bar{y}))|^2\dif
s\right)^{\frac{p}{2}}. \de Noting that \ce
|F(\bar{Y}_t)|\leq C(1+|\bar{Y}_t|),\\
|G(\bar{Y}_t)|\leq C(1+|\bar{Y}_t|),\\
|J(\bar{Y}_t,u)|\leq C(1+|\bar{Y}_t|)u, \de we obtain that \ce
\mE\left(\sup\limits_{0\leq t\leq T}|\bar{Y}_t|^p\right)
&\leq&4^{p-1}|\bar{y}|^p+4^{p-1}C\mE\left(\int^T_0(1+|\bar{Y}_s|^p)\dif s\right)\\
&\leq&4^{p-1}C(1+|\bar{y}|^p)+4^{p-1}C\mE\left(\int^T_0\left(\sup\limits_{0\leq
s\leq t}|\bar{Y}_s|^p\right)\dif t\right). \de
Thus, Gronwall's
inequality leads us to
 \ce \mE\left(\sup\limits_{0\leq
t\leq
T}|\bar{Y}_t|^p\right)\leq4^{p-1}C(1+|\bar{y}|^p)e^{4^{p-1}CT}.
\de

For $\bar{y}^1=\left(\begin{array}{c} y^1\\\eta^1 \end{array}
\right)$ and $\bar{y}^2=\left(\begin{array}{c} y^2\\\eta^2
\end{array} \right)$, by the similar deduction to the above one, we
can have that
\ce \mE\left(\sup\limits_{0\leq t\leq
T}|\bar{Y}^1_t-\bar{Y}^2_t|^p\right)\leq C|\bar{y}^1-\bar{y}^2|^p.
\de
By the Kolmogorov's continuity criterion in \cite{fk}, it yields
that for almost all $\omega\in\Omega$, $\bar{y}\mapsto
\bar{Y}_t(\bar{y})$ is continuous on $\mR^d\times\mR$ for all $t\geq
0$. This completes the proof.
\end{proof}

Take a two-sided $m$-dimensional L\'evy process \ce
\hat{L}_t=\left\{\begin{array}{l}
bt+AW_t+\int^t_0\int_{|u|\leq\delta}u\,\tilde{N}_\kappa(\dif s,\dif u), ~~\qquad\qquad t\geq0,\\
-b^-t-A^-W^-_t-\int_t^0\int_{|u|\leq\delta}u\tilde{N}_{\kappa^-}(\dif
s, \dif u), ~\quad t<0.
\end{array}
\right. \de
We consider the following stochastic integral
equations: for $x\in\mR^d, \tau\in\mR$, \be
&&h_t^{x,\tau}=x+e^{-\tau}\int_{-\infty}^te^s\sigma_i(h_s^{x,\tau})\dif \hat{L}^i_s, \quad\qquad t\in\mR,\label{au}\\
&&\tilde{h}_t^{x,\tau}=x+e^{-\tau}\int_0^t\sigma_i(\tilde{h}_s^{x,\tau})\dif
\tilde{\hat{L}}^i_s, \quad\qquad\qquad t>0, \label{trau} \ee where
$\tilde{\hat{L}}_t:=\int_{-\infty}^{\frac{1}{2}\log2t}e^s\dif
\hat{L}_s$. So,
$\tilde{h}_t^{x,\tau}=h_{\frac{1}{2}\log2t}^{x,\tau}.$

By Definition \ref{levypro1}, we know that $\tilde{\hat{L}}_t$ is a
L\'evy process. Thus, Proposition \ref{diffeom} implies
that $x\mapsto\tilde{h}_t^{x,\tau}$ is an a.s. $\cC^1$-diffeomorphism of $\mR^d$.   Then $x\mapsto h_t^{x,\tau}$ is also
a diffeomorphism of $\mR^d$ and $h_t^{x,\tau}$ is differentiable in
$\tau$. Set \ce H_t(\omega,x):=h_t^{x,\tau}|_{\tau=t},
\qquad\Gamma_t(\omega,x):=\frac{\partial h_t^{x,\tau}}{\partial
\tau}|_{\tau=t}, \qquad t\in\mR. \de

\bp\label{ou} The following results hold:

(i) \ce \dif H_t=\Gamma_t\dif t+\sigma_i(H_t)\dif \hat{L}^i_t, \de

(ii) \ce
H_{s+t}(\omega,x)=H_t(\theta_s\omega,x), \qquad
\Gamma_{s+t}(\omega,x)=\Gamma_t(\theta_s\omega,x), \qquad s,
t\in\mR,
\de
for a.s. $\omega$.
\ep
\begin{proof}
Set
$$
D(t,\tau):=\frac{\partial h_t^{x,\tau}}{\partial \tau}.
$$
Then $D(t,t)=\Gamma_t$, and it satisfies the following equation
\ce
D(t,\tau)&=&-e^{-\tau}\int_{-\infty}^te^s\sigma_i(h_s^{x,\tau})\dif \hat{L}^i_s\\
&&+e^{-\tau}\int_{-\infty}^te^s\frac{\partial\sigma_i}{\partial
x}(h_s^{x,\tau})D(s,\tau)\dif \hat{L}^i_s. \de So, for $s\leq t$,
$s, t\in\mR$, \ce
H_t&=&x+e^{-t}\int_{-\infty}^te^r\sigma_i(h_r^{x,t})\dif \hat{L}^i_r\\
&=&x+e^{-t}\int_{-\infty}^se^r\sigma_i(h_r^{x,t})\dif \hat{L}^i_r+e^{-t}\int_s^te^r\sigma_i(h_r^{x,t})\dif \hat{L}^i_r\\
&=&x+e^{-t}\int_{-\infty}^se^r\sigma_i(h_r^{x,t})\dif \hat{L}^i_r+\int_s^t\sigma_i(h_r^{x,r})\dif \hat{L}^i_r\\
&&-\int_s^t\left(\int_s^ue^{r-u}\sigma_i(h_r^{x,u})\dif \hat{L}^i_r\right)\dif u\\
&&+\int_s^t\left(\int_s^ue^{r-u}\frac{\partial\sigma_i}{\partial x}(h_r^{x,u})D(r,u)\dif \hat{L}^i_r\right)\dif u\\
&=&x+\int_s^t\Gamma_u\dif u+\int_s^t\sigma_i(H_r)\dif \hat{L}^i_r+e^{-t}\int_{-\infty}^se^r\sigma_i(h_r^{x,t})\dif \hat{L}^i_r\\
&&+\int_s^te^{-u}\left(\int_{-\infty}^se^r\sigma_i(h_r^{x,u})\dif \hat{L}^i_r\right)\dif u\\
&&-\int_s^te^{-u}\left(\int_{-\infty}^se^r\frac{\partial\sigma_i}{\partial
x}(h_r^{x,u})D(r,u)\dif \hat{L}^i_r\right)\dif u,
\de
where we have used the following formula
\be e^{-t}\int_s^te^r\sigma_i(h_r^{x,t})\dif
\hat{L}^i_r-\int_s^t\sigma_i(h_r^{x,r})\dif \hat{L}^i_r
&=&\int_s^t\left(\int_s^ue^{r-u}\frac{\partial\sigma_i}{\partial x}(h_r^{x,u})D(r,u)\dif \hat{L}^i_r\right)\dif u\no\\
&&-\int_s^t\left(\int_s^ue^{r-u}\sigma_i(h_r^{x,u})\dif
\hat{L}^i_r\right)\dif u. \label{form}
\ee
The proof of this formula   is   in Appendix.

Via integration by parts,
 \ce
H_t&=&x+\int_s^t\Gamma_u\dif u+\int_s^t\sigma_i(H_r)\dif \hat{L}^i_r+e^{-s}\int_{-\infty}^se^r\sigma_i(h_r^{x,s})\dif \hat{L}^i_r\\
&=&H_s+\int_s^t\Gamma_u\dif u+\int_s^t\sigma_i(H_r)\dif \hat{L}^i_r.
\de
Thus, \ce \dif H_t=\Gamma_t\dif t+\sigma_i(H_t)\dif\hat{L}^i_t.
\de

To (ii), for $s,t\in\mR$ and $\omega\in\Omega$,
\ce
H_{s+t}(\omega,x)&=&x+e^{-(s+t)}\int_{-\infty}^{s+t}e^r\sigma_i(H_r(\omega,x))\dif\hat{L}^i_r\\
&=&x+e^{-(s+t)}\int_{-\infty}^te^{s+u}\sigma_i(H_{s+u}(\omega,x))\dif(\hat{L}^i_{s+u}-\hat{L}^i_s)\\
&=&x+e^{-t}\int_{-\infty}^te^{u}\sigma_i(H_{s+u}(\omega,x))\dif
\hat{L}^i_{u}(\theta_s\omega). \de By the uniqueness of the solution
to Eq.(\ref{au}),   we get
\ce
H_{s+t}(\omega,x)=H_t(\theta_s\omega,x)
\de
for a.s. $\omega$ and $s,t\in\mR$. By a similar perfection procedure to one in \cite{ks},
we obtain the first result of (ii). The similar deduction to the above one admits us to have that
$\Gamma_{s+t}(\omega,x)=\Gamma_t(\theta_s\omega,x)$.
\end{proof}

By the above Proposition, we know that $H_t$ and $\Gamma_t$ are
stationary processes.

Now we consider the SDE \be \dif X_t=a(X_t)\dif
t+\sigma_i(X_t)\dif \hat{L}^i_t, \quad t\in\mR, \label{eq1} \ee and
the RDE \be \dif Y_t=\left(\frac{\partial}{\partial
x}H_t\right)^{-1}(\omega,Y_t)\[a(H_t(\omega,Y_t))-\Gamma_t(\omega,Y_t)\]\dif
t, \quad t\in\mR. \label{eq2} \ee

\bp\label{conjsde} Assume that $a(x)\in\cC_b^{1,\gamma}$,
$\sigma_i(x)\in\cC_b^{1,\gamma}$ for $i=1,2,\cdots,m$, the mapping
\ce \left(\begin{array}{c} y\\\eta \end{array}
\right)\mapsto\left(\begin{array}{c} y\\\eta \end{array}
\right)+\left(\begin{array}{c} \eta\sigma_i(y)u^i\\0 \end{array}
\right) \de is homeomorphic for $y\in\mR^d, \eta\in\mR,
u\in\{u\in\mR^m: |u|\leq\delta\}$ and the Jacobian matrix \ce
I+\left(\begin{array}{c} \eta\frac{\partial(\sigma_i(y))^1}{\partial
y_1}u^i
\quad \eta\frac{\partial(\sigma_i(y))^1}{\partial y_2}u^i ~\cdots ~ \eta\frac{\partial(\sigma_i(y))^1}{\partial y_d}u^i \quad (\sigma_i(y))^1u^i\\
\eta\frac{\partial(\sigma_i(y))^2}{\partial y_1}u^i
\quad \eta\frac{\partial(\sigma_i(y))^2}{\partial y_2}u^i ~\cdots ~ \eta\frac{\partial(\sigma_i(y))^2}{\partial y_d}u^i \quad (\sigma_i(y))^2u^i\\
\vdots ~\qquad\qquad \vdots \qquad \cdots \qquad \vdots \qquad\qquad \vdots\\
\eta\frac{\partial(\sigma_i(y))^d}{\partial y_1}u^i
\quad \eta\frac{\partial(\sigma_i(y))^d}{\partial y_2}u^i ~\cdots ~ \eta\frac{\partial(\sigma_i(y))^d}{\partial y_d}u^i \quad (\sigma_i(y))^du^i\\
0 ~\qquad\qquad 0 \qquad \cdots \qquad 0 \qquad\qquad 0
\end{array}\right)
\de is invertible for any $y,\eta$ a.e. and
$u\in\{u\in\mR^m:|u|\leq\delta\}$. Then

(i) the solution of Eq.(\ref{eq1}) generates a global cocycle
$\varphi_t$,

(ii) the solution of Eq.(\ref{eq2}) generates a global cocycle
$\psi_t$,

(iii) RDSs $\varphi_t$ and $\psi_t$ are conjugate with the
cohomology $H_0$. \ep
\begin{proof}
By Theorem 3.1 in \cite{kuni}, the solution of Eq.(\ref{eq1})
generates a global cocycle $\varphi_t$. However, by the assumptions
of $a(x)$ and $\sigma_i(x)$, we only obtain that Eq.(\ref{eq2})
generates a local cocycle denoted by $\psi_t$. To prove that
$\psi_t$ is global, we use conjugacy.

By the It\^o-Ventzell formula in Appendix for $t\geq 0$, we
have that \ce \dif H_t(\omega,Y_t)&=&\left(\frac{\partial}{\partial
x}H_t\right)\left(\frac{\partial}{\partial
x}H_t\right)^{-1}(\omega,Y_t)\[a(H_t(\omega,Y_t))-\Gamma_t(\omega,Y_t)\]\dif t\\
&&\Gamma_t(\omega,Y_t)\dif t+\sigma_i(H_t(\omega,Y_t))b^i\dif t+\sigma_i(H_t(\omega,Y_t))A^{ij}\dif W_t^j\\
&&+\int_{|u|\leq\delta}\sigma_i(H_t(\omega,Y_t))u^i\,\tilde{N}_\kappa(\dif t,\dif u)\\
&=&a(H_t(\omega,Y_t))\dif t+\sigma_i(H_t(\omega,Y_t))\dif
\hat{L}^i_t.
\de
For $t<0$, based on the similar deduction to the above one, the same result
holds. Thus,
\ce
\varphi_t(\omega,\cdot)=H_0(\theta_t\omega,\cdot)\circ\psi_t(\omega,\cdot)\circ
H_0(\omega,\cdot)^{-1}, \de and \ce
\psi_t(\omega,\cdot)=H_0(\theta_t\omega,\cdot)^{-1}\circ\varphi_t(\omega,\cdot)\circ
H_0(\omega,\cdot).
\de
Since $\varphi_t(\omega,\cdot)$ is global, so
is $\psi_t(\omega,\cdot)$.
\end{proof}

\subsection{Linearization of SDEs}\label{linear}

In this subsection, we study the relation between a SDE  and its  linearized SDE. A stochastic
version of the Hartman-Grobman theorem  serves our purpose (\cite{tmr, il2}). Using the promising technique of
conjugacy for flows generated by Stratonovitch SDEs with Brownian
motions and flows generated by RDEs, Imkeller and Lederer in
\cite{il2} proved a Hartman-Grobman theorem for continuous cocycles. Here
 we prove a
Hartman-Grobman theorem for discontinuous cocycles.

\bt\label{normsde} Assume that $a(x)\in\cC_b^{2,\gamma}$,
$\sigma_i(x)\in\cC_b^{2,\gamma}$ for $i=1,2,\cdots,m$, the mapping
\ce \left(\begin{array}{c} y\\\eta \end{array}
\right)\mapsto\left(\begin{array}{c} y\\\eta \end{array}
\right)+\left(\begin{array}{c} \eta\sigma_i(y)u^i\\0 \end{array}
\right) \de is homeomorphic for $y\in\mR^d, \eta\in\mR,
u\in\{u\in\mR^m: |u|\leq\delta\}$ and the Jacobian matrix \ce
I+\left(\begin{array}{c} \eta\frac{\partial(\sigma_i(y))^1}{\partial
y_1}u^i
\quad \eta\frac{\partial(\sigma_i(y))^1}{\partial y_2}u^i ~\cdots ~ \eta\frac{\partial(\sigma_i(y))^1}{\partial y_d}u^i \quad (\sigma_i(y))^1u^i\\
\eta\frac{\partial(\sigma_i(y))^2}{\partial y_1}u^i
\quad \eta\frac{\partial(\sigma_i(y))^2}{\partial y_2}u^i ~\cdots ~ \eta\frac{\partial(\sigma_i(y))^2}{\partial y_d}u^i \quad (\sigma_i(y))^2u^i\\
\vdots ~\qquad\qquad \vdots \qquad \cdots \qquad \vdots \qquad\qquad \vdots\\
\eta\frac{\partial(\sigma_i(y))^d}{\partial y_1}u^i
\quad \eta\frac{\partial(\sigma_i(y))^d}{\partial y_2}u^i ~\cdots ~ \eta\frac{\partial(\sigma_i(y))^d}{\partial y_d}u^i \quad (\sigma_i(y))^du^i\\
0 ~\qquad\qquad 0 \qquad \cdots \qquad 0 \qquad\qquad 0
\end{array}
\right) \de is invertible for any $y,\eta$ a.e. and
$u\in\{u\in\mR^m: |u|\leq\delta\}$. Suppose $a(0)=0$,
$\sigma_i(0)=0$ for $i=1,2,\cdots,m$, and set
$B_0:=\frac{\partial}{\partial x}a(0)$,
$B_i:=\frac{\partial}{\partial x}\sigma_i(0)$ for $i=1,2,\cdots,m$.

$\varphi_t(\omega,\cdot)$ still denotes the RDS generated by
Eq.(\ref{eq1}). Let $\varphi^0_t(\omega,\cdot)$ be a RDS generated
by \be \dif x_t=B_0x_t\dif t+B_ix_t\dif \hat{L}^i_t, \qquad t\in\mR
\label{lsde} \ee and suppose that all Lyapunov exponents of
$\varphi^0_t(\omega,\cdot)$ are non-zero. Then there exists a
measurable mapping $\varsigma: \Omega\times\mR^d\mapsto\mR^d$
satisfying the following properties

(i) $x\mapsto\varsigma(\omega,x)$ is a homeomorphism of $\mR^d$ and
$\varsigma(\omega,0)=0$,

(ii) for $\omega\in\Omega$, it holds that \ce
\varphi_t(\omega,x)=\varsigma(\theta_t\omega,\cdot)\circ
\varphi^0_t(\omega,\cdot)\circ\varsigma(\omega,x)^{-1}, \qquad
t\in\mR.
 \de
\et
\begin{proof}
{\bf Step 1.} Set
$$
F(\omega,y):=\left(\frac{\partial}{\partial
x}H_0\right)^{-1}(\omega,y)\[a(H_0(\omega,y))-\Gamma_0(\omega,y)\].
$$
Then Eq.(\ref{eq2}) can be rewritten as \ce
\quad\qquad\qquad\qquad\qquad\qquad\qquad\left\{\begin{array}{l}
\dif Y_t=F(\theta_t\omega,Y_t)\dif t, \quad t\in\mR, \\
Y_0=y.
\end{array}
\right. \qquad\qquad\qquad\qquad\qquad \eqref{eq2}' \de Moreover,
$F(\cdot,0)=0$ and for $\omega\in\Omega$,
$F(\omega,\cdot)\in\cC^1(\mR^d)$ by the assumptions for $a(x)$ and
$\sigma_i(x)$, $i=1,2,\cdots,m$.

Define
$$
f(\omega):=\left(\frac{\partial}{\partial y}F\right)(\omega,0).
$$
Hence,
$$
f(\omega)=\left(\frac{\partial}{\partial
x}H_0\right)^{-1}(\omega,0)\left[B_0\left(\frac{\partial}{\partial
x}H_0\right)(\omega,0) -\left(\frac{\partial}{\partial
x}\Gamma_0\right)(\omega,0)\right].
$$
Introduce the following linear equation: \be\left\{\begin{array}{l}
\dif y_t=f(\theta_t\omega)y_t\dif t, \quad t\in\mR,\\
y_0=y. \label{leq}
\end{array}
\right. \ee By the similar result to Theorem 3.2 in \cite{il2} (Only
the continuity of $\theta_t\omega$ in $t$ is changed into right
continuity with left limit), there exist measurable mappings
$\varrho:\Omega\mapsto(0,\infty)$ and
$\zeta:\Omega\times\mR^d\mapsto\mR^d$ such that

(i) $y\mapsto\zeta(\omega,y)$ is a homeomorphism of $\mR^d$ and
$\zeta(\omega,0)=0$,

(ii) for $\omega\in\Omega$, it holds that \ce
\psi_t(\omega,y)=\zeta(\theta_t\omega,\cdot)\circ
\psi^0_t(\omega,\cdot)\circ\zeta(\omega,y)^{-1}, \quad
\tau^1_-(\omega,y)\leq t\leq \tau^1_+(\omega,y), \de where
$\psi_t(\omega,\cdot)$ and $\psi^0_t(\omega,\cdot)$ are the RDS
generated by (\ref{eq2})$^{\prime}$ and (\ref{leq}), respectively,
and \ce
&&\tau^1_-(\omega,y):=\inf\{t<0:|\psi_s(\omega,y)|\leq\varrho(\theta_s\omega)~\mbox{for all}~t\leq s\leq 0\},\\
&&\tau^1_+(\omega,y):=\sup\{t>0:|\psi_s(\omega,y)|\leq\varrho(\theta_s\omega)~\mbox{for
all}~0\leq s\leq t\}. \de

Next, for $t\geq\tau^1_+(\omega,y)$, take
$\psi_{\tau^1_+}(\omega,y)$ and $\psi^0_{\tau^1_+}(\omega,y)$ as the
initial values of (\ref{eq2}) and (\ref{leq}), respectively, and
consider (\ref{eq2}) and (\ref{leq}) again. Using the cocycle
property of solutions to (\ref{eq2}) and (\ref{leq}), we obtain that
\ce
\psi_t(\omega,y)&=&\psi_{t-\tau^1_+}(\theta_{\tau^1_+}\omega,\psi_{\tau^1_+}(\omega,y))\\
&=&\zeta(\theta_{t-\tau^1_+}\theta_{\tau^1_+}\omega,\cdot)\circ
\psi^0_{t-\tau^1_+}(\theta_{\tau^1_+}\omega,\cdot)\circ\zeta(\theta_{\tau^1_+}\omega,\psi_{\tau^1_+}(\omega,y))^{-1}\\
&=&\zeta(\theta_t\omega,\cdot)\circ
\psi^0_{t-\tau^1_+}(\theta_{\tau^1_+}\omega,\cdot)\circ\zeta(\theta_{\tau^1_+}\omega,\cdot)^{-1}\circ\zeta(\theta_{\tau^1_+}\omega,\cdot)\\
&&\circ\psi^0_{\tau^1_+}(\omega,\cdot)\circ\zeta(\omega,y)^{-1}\\
&=&\zeta(\theta_t\omega,\cdot)\circ\psi^0_{t-\tau^1_+}(\theta_{\tau^1_+}\omega,\cdot)\circ
\psi^0_{\tau^1_+}(\omega,\cdot)\circ\zeta(\omega,y)^{-1}\\
&=&\zeta(\theta_t\omega,\cdot)\circ\psi^0_t(\omega,\cdot)\circ\zeta(\omega,y)^{-1},
\de for $\tau^1_+(\omega,y)\leq t\leq \tau^2_+(\omega,y)$, where \ce
\tau^2_+(\omega,y)=\sup\{t>\tau^1_+:|\psi_s(\omega,y)|\leq\varrho(\theta_s\omega)~\mbox{for
all}~\tau^1_+\leq s\leq t\}.
\de
Doing   this for
$t\geq\tau^2_+(\omega,y)$ and $t\leq\tau^1_-(\omega,y)$, we
 see that for $\omega\in\Omega$, \ce
\psi_t(\omega,y)=\zeta(\theta_t\omega,\cdot)\circ
\psi^0_t(\omega,\cdot)\circ\zeta(\omega,y)^{-1}, \qquad t\in\mR. \de

{\bf Step 2.} Consider the cocycles $\varphi^0_t(\omega,\cdot)$ and
$\psi^0_t(\omega,\cdot)$ generated by Eq.(\ref{lsde}) and
Eq.(\ref{leq}), respectively. By the It\^o-Ventzell formula in
Appendix for $t\geq0$, we get that \ce \dif
\left(\frac{\partial}{\partial
x}H_t\right)(\cdot,0)y_t&=&\left(\frac{\partial}{\partial
x}H_t\right)(\cdot,0)
\left(\frac{\partial}{\partial x}H_t\right)^{-1}(\cdot,0)\[B_0\left(\frac{\partial}{\partial x}H_t\right)(\cdot,0)\\
&&-\left(\frac{\partial}{\partial
x}\Gamma_t\right)(\cdot,0)\]y_t\dif t
+\left(\frac{\partial}{\partial x}\Gamma_t\right)(\cdot,0)y_t\dif t\\
&&+B_i\left(\frac{\partial}{\partial x}H_t\right)(\cdot,0)y_t\dif \hat{L}^i_t\\
&=&B_0\left(\frac{\partial}{\partial x}H_t\right)(\cdot,0)y_t\dif
t+B_i\left(\frac{\partial}{\partial x}H_t\right)(\cdot,0)y_t\dif
\hat{L}^i_t.
\de
For $t<0$, in terms of the similar deduction to the above one, the result also
holds.

Thus   $\varphi^0_t(\omega,\cdot)$ and $\psi^0_t(\omega,\cdot)$
are conjugate. $\left(\frac{\partial}{\partial
x}H_0\right)^{-1}(\cdot,0)$ is a cohomology of
$\psi^0_t(\omega,\cdot)$ and $\varphi^0_t(\omega,\cdot)$.

{\bf Step 3.} Combining Proposition \ref{conjsde} with {\bf Steps
1-2}, we obtain that \ce
\varphi_t(\omega,\cdot)&=&H_0(\theta_t\omega,\cdot)\circ\psi_t(\omega,\cdot)H_0(\omega,\cdot)^{-1}\\
&=&H_0(\theta_t\omega,\cdot)\circ\zeta(\theta_t\omega,\cdot)\circ\psi^0_t(\omega,\cdot)\circ\zeta(\omega,\cdot)^{-1}\circ H_0(\omega,\cdot)^{-1}\\
&=&H_0(\theta_t\omega,\cdot)\circ\zeta(\theta_t\omega,\cdot)\circ\left(\frac{\partial}{\partial
x}H_0\right)^{-1}(\theta_t\omega,0)\circ
\varphi^0_t(\omega,\cdot)\\
&&\circ\left(\frac{\partial}{\partial
x}H_0\right)(\omega,0)\circ\zeta(\omega,\cdot)^{-1}\circ
H_0(\omega,\cdot)^{-1}. \de
Set
$$
\varsigma(\omega,\cdot):=H_0(\omega,\cdot)\circ\zeta(\omega,\cdot)\circ\left(\frac{\partial}{\partial
x}H_0\right)^{-1}(\omega,0)\cdot.
$$
We thus have
 \ce
\varphi_t(\omega,x)=\varsigma(\theta_t\omega,\cdot)\circ
\varphi^0_t(\omega,\cdot)\circ\varsigma(\omega,x)^{-1}.
\de
The proof is completed.
\end{proof}

\br\label{local} Since $\psi_t$ and $\psi^0_t$ are locally conjugate
via $\zeta$ in {\bf Step 1}, the result in the above theorem is
local.
\er

We present an example to demonstrate the above theorem.

\bx\label{lex}
Take $a(x)=\beta x-x^l$ and $\sigma(x)=x$, where $\beta=\alpha+\sigma^2/2$. Here $\alpha$ is a real parameter and
$\sigma$ is a positive parameter. Let $l>1$ be an integer. Consider the following scalar nonlinear stochastic equation
\be
\dif X_t=(\beta X_t-X_t^l)\dif t+X_t\dif \hat{L}_t, \quad t\in\mR.
\label{exeq1}
\ee
Lyapunov stability for this equation has been studied in a special case where the L\'evy process is replaced by a compound Poisson process in \cite{g}.

In this example, $B_0=\frac{\partial}{\partial x}a(0)=\beta$, and $B=\frac{\partial}{\partial x}\sigma(0)=1$.
Thus, by Theorem \ref{normsde}, the RDS $\varphi_t(\omega,\cdot)$ generated by Eq.(\ref{exeq1})
is conjugate to the RDS $\varphi^0_t(\omega,\cdot)$ generated by the linear stochastic equation
\be
\dif x_t=\beta x_t\dif t+x_t\dif \hat{L}_t, \quad t\in\mR.
\label{exeq2}
\ee


\ex

\section{Conjugacy for Marcus SDEs and RDEs, and existence of global attractors}\label{marcussde}

In this section, we first introduce Marcus SDEs, then we prove the
conjugacy for a Marcus SDE and a related RDE  and the    existence
of the global attractor for the Marcus SDE.

\subsection{Conjugacy for Marcus SDEs and RDEs}\label{conmarsde}

A Marcus canonical SDE was introduced in \cite{ma}. For the
L\'evy process $\hat{L}_t$  in Section \ref{conitosde}, this SDE is
  as follows:
  \be \dif
\bar{X}_t=\bar{a}(\bar{X}_t)\dif
t+\bar{\sigma}_i(\bar{X}_t)\diamond\dif \hat{L}^i_t, \qquad t\in\mR,
\label{marsde} \ee where $\bar{a}, \bar{\sigma}_1, \bar{\sigma}_2,
\cdots, \bar{\sigma}_m$ are vector fields on $\mR^d$, and
``$\diamond$" denotes the Marcus differential (\cite{da,kuni}).
The precise definition is as follows: \ce \dif
\bar{X}_t=\left\{\begin{array}{l}
\bar{a}(\bar{X}_t)\dif t+\bar{\sigma}_i(\bar{X}_t)\circ\dif \hat{L}^i_c(t)+\bar{\sigma}_i(\bar{X}_{t-})\dif \hat{L}^i_d(t)\\
+\Big\{\Phi(\Delta \hat{L}_t,\bar{X}_{t-})-\bar{X}_{t-}-\bar{\sigma}_i(\bar{X}_{t-})\Delta \hat{L}^i_t\Big\}, \quad t\geq0,\\
\bar{a}(\bar{X}_t)\dif t+\bar{\sigma}_i(\bar{X}_t)\circ\dif \hat{L}^i_c(t)+\bar{\sigma}_i(\bar{X}_{t-})\dif \hat{L}^i_d(t)\\
+\Big\{\Phi(\Delta
\hat{L}_t,\bar{X}_{t+})-\bar{X}_{t+}-\bar{\sigma}_i(\bar{X}_{t+})\Delta
\hat{L}^i_t\Big\}, \quad t\leq0,
\end{array}
\right. \de where ``$\circ$" stands for the Stratonovitch
differential, $\hat{L}_c(t)$ and $\hat{L}_d(t)$ are continuous and
purely discontinuous parts of the L\'evy process $\hat{L}(t)$,
respectively, \ce \Delta \hat{L}_t:=\left\{\begin{array}{l}
\hat{L}_t-\hat{L}_{t-}, \quad t\geq0,\\
\hat{L}_t-\hat{L}_{t+}, \quad t\leq0,
\end{array}
\right. \de and $\Phi(z,x)$ is the solution flow of the following
ordinary differential equation \be\left\{\begin{array}{l}
\frac{\partial\Phi(z,x)}{\partial z_i}=\bar{\sigma}_i(\Phi(z,x)), \qquad i=1,2,\cdots,m,\\
\Phi(0,x)=x.
\end{array}
\right. \label{pardif} \ee

\br\label{flowde} If the L\'evy process $\hat{L}_t$ has no
discontinuous part, Eq.(\ref{marsde}) changes into a Stratonovitch
SDE, and $\Phi(z,x)$ is the same as the auxiliary function in the
Doss-Sussmann flow decomposition (\cite{il2,is}). \er

Assume that $\bar{\sigma}_i(x)$ is smooth for $i=1,2,\cdots,m$
and the Lie bracket $[\bar{\sigma}_i,\bar{\sigma}_j] =0$ in the sense of
differential geometry for $i,j=1,2,\cdots,m, i\neq j$. Then
Eq.(\ref{pardif}) possesses a solution denoted by $\Phi(z,x)$ and
$x\mapsto\Phi(z,x)$ is a diffeomorphism of $\mR^d$, as in \cite{sm}.

For $\mu>0$, introduce the following process
$$
Z_t=e^{-\mu t}\int_{-\infty}^te^{\mu s}\dif \hat{L}_s, \quad
t\in\mR.
$$

\bl\label{ou} $Z_t$ has the following properties:

(i) \ce \dif Z_t=-\mu Z_t\dif t+\dif \hat{L}_t, \quad t\in\mR, \de

(ii) \ce Z_{s+t}(\omega)=Z_t(\theta_s\omega), \quad s,t\in\mR. \de
\el
\begin{proof}
Set $\bar{\hat{L}}_t:=\int_{-\infty}^te^{\mu s}\dif \hat{L}_s$, and
then $Z_t=e^{-\mu t}\bar{\hat{L}}_t$. By the It\^o-Ventzell formula
in Appendix for $t\geq0$, one obtains that \ce
\dif Z_t&=&e^{-\mu t}e^{\mu t}b\dif t+e^{-\mu t}e^{\mu t}A\dif W_t-\mu e^{-\mu t}\bar{\hat{L}}_t\dif t\\
&&+\int_{|u|\leq\delta}\left(e^{-\mu t}(\bar{\hat{L}}_t+e^{\mu t}u)-e^{-\mu t}\bar{\hat{L}}_t\right)\,\tilde{N}_\kappa(\dif t,\dif u)\\
&&+\int_{|u|\leq\delta}\left(e^{-\mu t}(\bar{\hat{L}}_t+e^{\mu t}u)-e^{-\mu t}\bar{\hat{L}}_t-e^{-\mu t}e^{\mu t}u\right)\nu(\dif u)\dif t\\
&=&-\mu Z_t\dif t+\dif \hat{L}_t. \de
For $t<0$, by the similar deduction to the above one, we could prove the same result. This
yields (i).

To prove (ii), it holds that
\ce
Z_{s+t}(\omega)&=&e^{-\mu(s+t)}\int_{-\infty}^{s+t}e^{\mu r}\dif \hat{L}_r\\
&=&e^{-\mu(s+t)}\int_{-\infty}^te^{\mu(s+u)}\dif (\hat{L}_{s+u}-\hat{L}_s)\\
&=&e^{-\mu t}\int_{-\infty}^te^{\mu u}\dif
\hat{L}_u(\theta_s\omega)=Z_t(\theta_s\omega). \de The proof is
completed.
\end{proof}

Denote
$$
\bar{H}_t(x):=\Phi(Z_t,x), \qquad t\in\mR.
$$
By the above Lemma,
$$
\bar{H}_{s+t}(\omega,x)=\bar{H}_t(\theta_s\omega,x)
$$
and by Theorem 4.4.7 or Theorem 4.4.28 in \cite{da}, \ce
\dif \bar{H}_t&=&\bar{\sigma}_i(\bar{H}_t)\diamond\dif Z^i_t\\
&=&-\mu\bar{\sigma}_i(\bar{H}_t)Z^i_t\dif
t+\bar{\sigma}_i(\bar{H}_t)\diamond\dif\hat{L}^i_t. \de Introduce
the following RDE \be \dif \bar{Y}_t=\left(\frac{\partial}{\partial
x}\bar{H}_t\right)^{-1}(\omega,\bar{Y}_t)\[\bar{a}\big(\bar{H}_t(\omega,\bar{Y}_t)\big)
+\mu\bar{\sigma}_i\big(\bar{H}_t(\omega,\bar{Y}_t)\big)Z_t^i\]\dif
t, \quad t\in\mR. \label{rdegakn} \ee

\bp\label{conjmarsde} RDSs $\bar{\varphi}_t(\omega,\cdot)$ and
$\bar{\psi}_t(\omega,\cdot)$ are conjugate with cohomology
$\bar{H}_0$, where $\bar{\varphi}_t(\omega,\cdot),
\bar{\psi}_t(\omega,\cdot)$ are RDSs generated by Eq.(\ref{marsde})
and Eq.(\ref{rdegakn}), respectively. \ep
\begin{proof}
We firstly show that the solutions of Eq.(\ref{marsde}) and
Eq.(\ref{rdegakn}) generate global RDSs. By Theorem 3.16 in
\cite{kuni}, the solution of Eq.(\ref{marsde}) generates a global
cocycle $\bar{\varphi}_t(\omega,\cdot)$. To prove that
$\bar{\psi}_t(\omega,\cdot)$ is global, we use conjugacy.

Rewrite $\dif \bar{H}_t$ as the It\^o form for $t\geq0$ \ce
\dif \bar{H}_t&=&-\mu\bar{\sigma}_i(\bar{H}_t)Z_t^i\dif t+\bar{\sigma}_i(\bar{H}_t)b^i\dif t+\bar{\sigma}_i(\bar{H}_t)A^{ij}\dif W_t^j\\
&&+\frac{1}{2}\left(\frac{\partial \bar{\sigma}_i}{\partial x_k}\right)(\bar{H}_t)\left(\bar{\sigma}_j(\bar{H}_t)\right)^kA^{il}A^{jl}\dif t\\
&&+\int_{|u|\leq\delta}\big(\Phi(u,\bar{H}_{t-})-\bar{H}_{t-}\big)\,\tilde{N}_\kappa(\dif t,\dif u)\\
&&+\int_{|u|\leq\delta}\big(\Phi(u,\bar{H}_{t-})-\bar{H}_{t-}-\bar{\sigma}_i(\bar{H}_{t-})u^i\big)\nu(\dif
u)\dif t, \de and for $t\leq0$ \ce
\dif \bar{H}_t&=&-\mu\bar{\sigma}_i(\bar{H}_t)Z_t^i\dif t-\bar{\sigma}_i(\bar{H}_t)(b^-)^i\dif t-\bar{\sigma}_i(\bar{H}_t)(A^-)^{ij}\dif (W^-)_t^j\\
&&+\frac{1}{2}\left(\frac{\partial \bar{\sigma}_i}{\partial x_k}\right)(\bar{H}_t)\left(\bar{\sigma}_j(\bar{H}_t)\right)^k(A^-)^{il}(A^-)^{jl}\dif t\\
&&+\int_{|u|\leq\delta}\big(\Phi(u,\bar{H}_{t+})-\bar{H}_{t+}\big)\,\tilde{N}_{\kappa^-}(\dif t,\dif u)\\
&&+\int_{|u|\leq\delta}\big(\Phi(u,\bar{H}_{t+})-\bar{H}_{t+}-\bar{\sigma}_i(\bar{H}_{t+})u^i\big)\nu^-(\dif
u)\dif t. \de Thus, the It\^o-Ventzell formula in Appendix for
$t\geq 0$ implies that
 \ce \dif
\bar{H}_t(\omega,\bar{Y}_t)&=&\left(\frac{\partial}{\partial
x}\bar{H}_t\right)\left(\frac{\partial}{\partial
x}\bar{H}_t\right)^{-1}(\omega,\bar{Y}_t)\[\bar{a}\big(\bar{H}_t(\omega,\bar{Y}_t)\big)+\mu\bar{\sigma}_i\big(\bar{H}_t(\omega,\bar{Y}_t)\big)Z_t^i\]\dif t\\
&&-\mu\bar{\sigma}_i(\bar{H}_t(\omega,Y_t))Z_t^i\dif
t+\bar{\sigma}_i(\bar{H}_t(\omega,\bar{Y}_t))b^i\dif t
+\bar{\sigma}_i(\bar{H}_t(\omega,\bar{Y}_t))A^{ij}\dif W_t^j\\
&&+\frac{1}{2}\left(\frac{\partial \bar{\sigma}_i}{\partial
x_k}\right)\big(\bar{H}_t(\omega,\bar{Y}_t)\big)\left(\bar{\sigma}_j(\bar{H}_t(\omega,\bar{Y}_t))\right)^k
A^{il}A^{jl}\dif t\\
&&+\int_{|u|\leq\delta}\left(\Phi(u,\bar{H}_{t-}(\omega,\bar{Y}_t))
-\bar{H}_{t-}(\omega,\bar{Y}_t)\right)\,\tilde{N}_\kappa(\dif t,\dif u)\\
&&+\int_{|u|\leq\delta}\left(\Phi(u,\bar{H}_{t-}(\omega,\bar{Y}_t))-\bar{H}_{t-}(\omega,\bar{Y}_t)
-\bar{\sigma}_i(\bar{H}_{t-}(\omega,\bar{Y}_t))u^i\right)\nu(\dif u)\dif t\\
&=&\bar{a}(\bar{H}_t(\omega,\bar{Y}_t))\dif
t+\bar{\sigma}_i(\bar{H}_t(\omega,\bar{Y}_t))\diamond\dif
\hat{L}^i_t. \de

For $t<0$, by the similar deduction to the above one, the same result
holds. Thus, \ce
\bar{\varphi}_t(\omega,\cdot)=\bar{H}_0(\theta_t\omega,\cdot)\circ\bar{\psi}_t(\omega,\cdot)\circ\bar{H}_0(\omega,\cdot)^{-1},
\de and \ce
\bar{\psi}_t(\omega,\cdot)=\bar{H}_0(\theta_t\omega,\cdot)^{-1}\circ\bar{\varphi}_t(\omega,\cdot)\circ
\bar{H}_0(\omega,\cdot). \de Since $\bar{\varphi}_t(\omega,\cdot)$
is global, so is $\bar{\psi}_t(\omega,\cdot)$.
\end{proof}

\medskip

We present two examples to explain what the cohomology
$\bar{H}_0$ is.

\bx\label{ex1}   Take $\bar{\sigma}_i(x)=\bar{\sigma}_ix$ for
$i=1,2,\cdots,m$, where $\bar{\sigma}_i\in\mR^{d\times d}$ and
$\bar{\sigma}_i\bar{\sigma}_j=\bar{\sigma}_j\bar{\sigma}_i$, $i\neq
j$.   Then Eq.(\ref{pardif}) changes into
\ce\left\{\begin{array}{l}
\frac{\partial\Phi(z,x)}{\partial z_i}=\bar{\sigma}_i\Phi(z,x), \qquad i=1,2,\cdots,m,\\
\Phi(0,x)=x.
\end{array}
\right. \de Solving it, we obtain that
$\Phi(z,x)=x\exp\{\bar{\sigma}_iz^i\}$. So, \ce
\bar{H}_t(\cdot,x)&=&\Phi(Z_t,x)=x\exp\{\bar{\sigma}_iZ_t^i\}=x\exp\{\bar{\sigma}_iZ_0^i\}\circ\theta_t\cdot\\
&=&\bar{H}_0(\theta_t\cdot,x). \de

In this case, the RDE
(\ref{rdegakn}) becomes \ce \dif
\bar{Y}_t=\exp\{-\bar{\sigma}_iZ_t^i\}\[\bar{a}\big(\bar{Y}_t\exp\{\bar{\sigma}_iZ_t^i\}\big)
+\mu\bar{\sigma}_i\bar{Y}_t\exp\{\bar{\sigma}_iZ_t^i\}Z_t^i\]\dif t,
\quad t\in\mR. \de \ex

\bx\label{ex2}   Take
$\bar{\sigma}_i(x)=\bar{\sigma}_ix+\beta_i$ for $i=1,2,\cdots,m$,
where $\beta_i\in\mR^{d}$, $\bar{\sigma}_i\in\mR^{d\times d}$ and
$\bar{\sigma}_i\bar{\sigma}_j=\bar{\sigma}_j\bar{\sigma}_i$, $\bar{\sigma}_i\beta_j=\bar{\sigma}_j\beta_i$, $i\neq
j$. Thus the solution of Eq.(\ref{pardif}) is \ce
\Phi(z,x)=x\exp\{\bar{\sigma}_iz^i\}+\bar{\sigma}^{-1}_i\big(\exp\{\bar{\sigma}_iz^i\}-I\big)\beta_i,
\de where $\bar{\sigma}^{-1}_i$ is the pseudo-inverse of the matrix
$\bar{\sigma}_i$.

Hence
 \ce
\bar{H}_t(\cdot,x)&=&\Phi(Z_t,x)=x\exp\{\bar{\sigma}_iZ_t^i\}+\bar{\sigma}^{-1}_i\big(\exp\{\bar{\sigma}_iZ_t^i\}-I\big)\beta_i\\
&=&\left[x\exp\{\bar{\sigma}_iZ_0^i\}+\bar{\sigma}^{-1}_i\big(\exp\{\bar{\sigma}_iZ_0^i\}-I\big)\beta_i\right]\circ\theta_t\cdot\\
&=&\bar{H}_0(\theta_t\cdot,x). \de

Finally, the RDE (\ref{rdegakn})
is written as \ce \dif
\bar{Y}_t=\exp\{-\bar{\sigma}_iZ_t^i\}\[\bar{a}\big(\bar{H}_t(\omega,\bar{Y}_t)\big)
+\mu\bar{\sigma}_i\bar{H}_t(\omega,\bar{Y}_t)Z_t^i+\mu\beta_iZ_t^i\]\dif
t, \quad t\in\mR. \de \ex

\subsection{Existence of global attractors}\label{exat}

We need   some   concepts. The function
$V:\mR^d\mapsto\mR_+$ is called the Lyapunov function of \be \dif
\bar{y}_t=\bar{a}(\bar{y}_t)\dif t, \label{deters} \ee if there
exists $\alpha>0$ such that
\be
\limsup\limits_{|y|\rightarrow\infty}\<\nabla\log
V(y),\bar{a}(y)\>\leq-\alpha. \label{lyacon} \ee A function
$U:\mR^d\mapsto\mR_+$ is said to {\it preserve temperedness} if
$U^{-1}(G)$ is tempered for tempered $G$. Finally, a function
$k:\mR^m\mapsto\mR$ is said to be {\it subexponentially growing} if
there exists $c>0$ such that $e^{-c|z|}k(z)$ is bounded.

\bp\label{rdegloatt} Let $V$ be a Lyapunov function of
Eq.(\ref{deters}) and assume that $V$ preserves temperedness. Set
\ce l(z,y):=\left(\frac{\partial \Phi}{\partial
y}\right)^{-1}(z,y)\bar{\sigma}_i\left(\Phi(z,y)\right)z^i, \quad
z\in\mR^m, y\in\mR^d. \de
Assume  that there are subexponentially
growing functions $k_1$ and $k_2$ such that
\be
\limsup\limits_{|y|\rightarrow\infty}\sup\limits_{z\in\mR^m}\left|\left<\nabla\log V(y),\frac{l(z,y)}{k_1(z)}\right>\right|=0,\label{c1}\\
\limsup\limits_{|y|\rightarrow\infty}\sup\limits_{z\in\mR^m}\frac{\left<\nabla\log
V(y),\bar{a}(y)
-(\frac{\partial \Phi}{\partial y})^{-1}(z,y)\bar{a}(\Phi(z,y))\right>}{\left|\left<\nabla\log V(y),\bar{a}(y)\right>\right|k_2(z)}\leq 1,\label{c2}
\ee
and
\ce
\lim\limits_{z\rightarrow0}k_2(z)=0.
\de
Then $\bar{\psi}_t$ has a
global random attractor. \ep

The proof is similar to one for Corollary 2.1 in \cite{is} and we thus
omit it.

Now, we state and prove the main result in this subsection.

\bt\label{sdegloatt} Under the conditions of Proposition
\ref{rdegloatt}, $\bar{\varphi}_t$ has a global random attractor.
\et
\begin{proof}
By the above Proposition, we know that $\bar{\psi}_t$ has a global
random attractor $A(\omega)$. Then by definition \ce
\bar{\psi}_t(\omega,A(\omega))=A(\theta_t\omega). \de
Define
$B(\omega):=\bar{H}_0(\omega,A(\omega))$.
Therefore
 \ce
\bar{\varphi}_t(\omega,B(\omega))&=&\bar{\varphi}_t(\omega,\cdot)\circ\bar{H}_0(\omega,A(\omega))\\
&=&\bar{\varphi}_t(\omega,\cdot)\circ\bar{H}_0(\omega,\cdot)\circ A(\omega)\\
&=&\bar{H}_0(\theta_t\omega,\cdot)\circ\bar{\psi}_t(\omega,\cdot)\circ\bar{H}_0(\omega,\cdot)^{-1}\circ\bar{H}_0(\omega,\cdot)\circ A(\omega)\\
&=&\bar{H}_0(\theta_t\omega,\cdot)\circ\bar{\psi}_t(\omega,\cdot)\circ A(\omega)\\
&=&\bar{H}_0(\theta_t\omega,A(\theta_t\omega))\\
&=&B(\theta_t\omega). \de

For any random tempered bounded set $G(\omega)$, by the assumption
on $\bar{\sigma}_i$ for $i=1,2,\cdots,m$,
$\bar{H}_0(\omega,\cdot)^{-1}\circ G(\omega)$ is tempered. Thus, by
Definition \ref{globattr}, \ce
\lim\limits_{t\rightarrow+\infty}dist(\bar{\psi}_t(\theta_{-t}\omega,\cdot)\circ\bar{H}_0(\theta_{-t}\omega,\cdot)^{-1}
\circ G(\theta_{-t}\omega),A(\omega))=0. \de Since, by conjugacy and
the assumption on $\bar{\sigma}_i$ for $i=1,2,\cdots,m$, \ce
&&dist(\bar{\varphi}_t(\theta_{-t}\omega,G(\theta_{-t}\omega)),B(\omega))\\
&=&dist(\bar{\varphi}_t(\theta_{-t}\omega,\cdot)\circ G(\theta_{-t}\omega),\bar{H}_0(\omega,A(\omega)))\\
&=&dist\left(\bar{H}_0(\omega,\cdot)\circ\bar{\psi}_t(\theta_{-t}\omega,\cdot)\circ\bar{H}_0(\theta_{-t}\omega,\cdot)^{-1}
\circ G(\theta_{-t}\omega),\bar{H}_0(\omega,A(\omega))\right)\\
&\leq&dist\left(\bar{\psi}_t(\theta_{-t}\omega,\cdot)\circ\bar{H}_0(\theta_{-t}\omega,\cdot)^{-1}
\circ G(\theta_{-t}\omega),A(\omega)\right), \de we obtain that \ce
\lim\limits_{t\rightarrow+\infty}dist(\bar{\varphi}_t(\theta_{-t}\omega,G(\theta_{-t}\omega)),B(\omega))=0,
\de for all tempered random bounded set $G(\omega)\subset\mR^d$.

Finally, by Definition \ref{globattr}, $B(\omega)$ is a global
random attractor for $\bar{\varphi}_t$. This completes the proof.
\end{proof}

 In the following example, we apply the above theorem to show the existence of a  global attractor  for the stochastic Duffing-van
der Pol equations with L\'evy processes.

\bx\label{ex3}
Consider the following Duffing-van der Pol system with L\'evy processes:
\ce\left\{\begin{array}{l}
\dif x_1(t)=x_2(t)\dif t,\\
\dif x_2(t)=[\gamma_1 x_1(t)+\gamma_2 x_2(t)-x_1(t)^3-x_1(t)^2x_2(t)]\dif t\\
\qquad\qquad +\sigma_1x_1(t)\diamond\dif \hat{L}^1_t+\sigma_2\diamond\dif \hat{L}^2_t,
\end{array}
\right. \de
where $\gamma_1, \gamma_2, \sigma_1, \sigma_2\in\mR^1$ and $\hat{L}_t=(\hat{L}^1_t,\hat{L}^2_t)$ is
a $2$-dimensional L\'evy process as defined in Section \ref{conitosde}. This system has been studied
in special cases: the L\'evy process is replaced by the Brownian motion in \cite{is,sh} and by the
Poisson process in \cite{x}. Set
\ce
&&\bar{X}_1(t):=x_1(t),\\
&&\bar{X}_2(t):=x_2(t)-\gamma_2 x_1(t)+\frac{1}{3}x_1(t)^3.
\de
Then the above system is transformed into
\ce
\dif\bar{X}_t=\bar{a}(\bar{X}_t)\dif
t+(\bar{\sigma}_i\bar{X}_t+\beta_i)\diamond\dif \hat{L}^i_t, \qquad t\in\mR,
\de
where
\ce
\bar{a}(y)=\left(\begin{array}{c} \gamma_2 y_1-\frac{1}{3}y_1^3+y_2 \\
\gamma_1 y_1-y_1^3 \end{array}\right),\quad\quad\\
\bar{\sigma}_1=\left(\begin{array}{c} 0 \quad 0 \\
\sigma_1 \quad 0\end{array}\right),
\quad \beta_1=0, \quad \bar{\sigma}_2=0, \quad \mbox{ and }  \beta_2=\left(\begin{array}{c} 0\\
\sigma_2 \end{array}\right).
\de

For the deterministic equation
\ce
\dif\bar{y}_t=\bar{a}(\bar{y}_t)\dif t,
\de
  the function
\ce
V(y)=\frac{7}{24}y_1^4+\frac{1}{4}y_1^2+\frac{1}{4}y_2^2+\frac{1}{2}(y_1-y_2)^2
\de
is a Lyapunov function. Indeed,
\ce
\nabla V(y)=\left(\begin{array}{c} \frac{7}{6}y_1^3+\frac{3}{2}y_1-y_2\\
\frac{3}{2}y_2-y_1 \end{array}\right),
\de
and
\ce
\<\nabla V(y),\bar{a}(y)\>&=&-\frac{7}{18}y_1^6+\left(\frac{7}{6}\gamma_2+\frac{1}{2}\right)y_1^4+\left(\frac{3}{2}\gamma_2-\gamma_1\right)y_1^2\\
&&+\left(\frac{3}{2}-\gamma_2+\frac{3}{2}\gamma_1\right)y_1y_2-y_2^2.
\de
Thus, there exists $\eta>0$ such that
\be
\limsup\limits_{|y|\rightarrow\infty}\frac{\<\nabla V(y),\bar{a}(y)\>}{\kappa(y)} \leq  -\eta,
\label{esti}
\ee
where $\kappa(y)=y_1^6+y_2^2$. Thus, $V(y)$ is a Lyapunov function.

Noticing that
\ce
\exp\{\bar{\sigma}_1z^1\}=\left(\begin{array}{c} 1 \qquad 0\\
\sigma_1z^1 \quad 1\end{array}\right), \quad
\exp\{-\bar{\sigma}_1z^1\}=\left(\begin{array}{c} 1 \quad\qquad 0\\
-\sigma_1z^1 \quad 1\end{array}\right),
\de
by Example \ref{ex2}, we get
\ce
\Phi(z,y)&=&y\exp\{\bar{\sigma}_1z^1\}+\beta_2z^2\\
&=&\left(\begin{array}{c} y_1\\
\sigma_1y_1z^1+y_2+\sigma_2z^2 \end{array}\right).
\de
Therefore,
\ce
l(z,y)&=&\left(\frac{\partial \Phi}{\partial y}\right)^{-1}(z,y)\left[\bar{\sigma}_1\left(\Phi(z,y)\right)z^1
+\bar{\sigma}_2\left(\Phi(z,y)\right)z^2\right]\\
&=&\exp\{-\bar{\sigma}_1z^1\}\left[\bar{\sigma}_1(y\exp\{\bar{\sigma}_1z^1\}+\beta_2z^2)z^1+\beta_2z^2\right]\\
&=&\bar{\sigma}_1yz^1+\exp\{-\bar{\sigma}_1z^1\}\beta_2z^2\\
&=&\left(\begin{array}{c} 0\\
\sigma_1y_1z^1+\sigma_2z^2 \end{array}\right),
\de
and
\ce
\<\nabla V(y),l(z,y)\>&=&\left(\frac{3}{2}y_2-y_1\right)\left(\sigma_1y_1z^1+\sigma_2z^2\right)\\
&=&\sigma_1\left(\frac{3}{2}y_1y_2-y^2_1\right)z^1+\sigma_2\left(\frac{3}{2}y_2-y_1\right)z^2.
\de
Set $k_1(z):=C|z|$. Then $k_1(z)$ is subexponentially growing and
\ce
\limsup\limits_{|y|\rightarrow\infty}\sup\limits_{z\in\mR^m}\left|\left<\nabla\log V(y),\frac{l(z,y)}{k_1(z)}\right>\right|=0.
\de
Thus, the condition (\ref{c1}) in Proposition \ref{rdegloatt} is satisfied.

To justify (\ref{c2}), we decompose $\bar{a}(y)$ into two parts:
\ce
\bar{a}(y)=Ay+\left(\begin{array}{c} -\frac{1}{3}y_1^3\\
-y_1^3 \end{array}\right),
\de
where $A:=\left(\begin{array}{c} \gamma_2 \quad 1\\
\gamma_1 \quad 0 \end{array}\right)$. Hence,
\ce
\bar{a}(y)-(\frac{\partial \Phi}{\partial y})^{-1}(z,y)\bar{a}(\Phi(z,y))
&=&Ay-(\frac{\partial \Phi}{\partial y})^{-1}(z,y)A\Phi(z,y)\\
&&+\left(\begin{array}{c} -\frac{1}{3}y_1^3\\
-y_1^3 \end{array}\right)-(\frac{\partial \Phi}{\partial y})^{-1}(z,y)\left(\begin{array}{c} -\frac{1}{3}y_1^3\\
-y_1^3 \end{array}\right)\\
&=&Ay-\exp\{-\bar{\sigma}_1z^1\}A\exp\{\bar{\sigma}_1z^1\}y-\exp\{-\bar{\sigma}_1z^1\}A\beta_2z^2\\
&&+(I-\exp\{-\bar{\sigma}_1z^1\})\left(\begin{array}{c} -\frac{1}{3}y_1^3\\
-y_1^3 \end{array}\right)\\
&=&\left(\begin{array}{c} -1 \qquad\qquad 0\\
\gamma_2+\sigma_1z^1 \quad 1\end{array}\right)y\sigma_1z^1-\left(\begin{array}{c} 1\\
-\sigma_1z^1\end{array}\right)\sigma_2z^2\\
&&+\left(\begin{array}{c} 0\\
-\frac{1}{3}y_1^3 \end{array}\right)\sigma_1z^1\\
&=:&S_1+S_2+S_3.
\de
Define
\ce
k_2(z):&=&(\sqrt{2}+|\gamma_2|+|\sigma_1||z|)|\sigma_1||z|\\
&&+(1+|\sigma_1||z|)|\sigma_2||z|+|\sigma_1||z|.
\de
Then $k_2(z)$ is subexponentially growing and
\ce
\lim\limits_{z\rightarrow0}k_2(z)=0.
\de
Moreover,
\ce
\frac{\left<\nabla\log V(y),\bar{a}(y)
-(\frac{\partial \Phi}{\partial y})^{-1}(z,y)\bar{a}(\Phi(z,y))\right>}{\left|\left<\nabla\log V(y),\bar{a}(y)\right>\right|k_2(z)}
&=&\frac{\left<\nabla\log V(y),S_1+S_2+S_3\right>}{\left|\left<\nabla\log V(y),\bar{a}(y)\right>\right|k_2(z)}\\
&\leq&\frac{|\nabla V(y)||y|}{\left|\left<\nabla V(y),\bar{a}(y)\right>\right|}
+\frac{|\nabla V(y)|}{\left|\left<\nabla V(y),\bar{a}(y)\right>\right|}\\
&&+\frac{|(\frac{3}{2}y_2-y_1)\frac{1}{3}y_1^3|}{\left|\left<\nabla V(y),\bar{a}(y)\right>\right|}.
\de
It follows from (\ref{esti}) that there exists a constant $C>0$ such that
\ce
\left|\left<\nabla V(y),\bar{a}(y)\right>\right|\leq C \; \kappa(y).
\de
Thus, by simple calculation, we obtain that
\ce
\limsup\limits_{|y|\rightarrow\infty}\sup\limits_{z\in\mR^m}\frac{\left<\nabla\log
V(y),\bar{a}(y)
-(\frac{\partial \Phi}{\partial y})^{-1}(z,y)\bar{a}(\Phi(z,y))\right>}{\left|\left<\nabla\log V(y),\bar{a}(y)\right>\right|k_2(z)}\leq 1.
\de
By Theorem \ref{sdegloatt}, the stochastic Duffing-van der Pol equation has a global attractor.
\ex

\section{Appendix}\label{itit}

In this section, we recall the It\^o-Ventzell formula in the context
of SDEs with L\'evy processes (see \cite{q1}), and then provide a
proof for the formula (\ref{form}) in Section \ref{itosde}.

\medskip

\textbf{The It\^o-Ventzell formula:}  \\
Let ($\mU,\cU$) be a measurable space and $n$ be a $\sigma$-finite
measure on it. Let $\mU_{0}$ be a set in $\cU$ such that
$n(\mU-\mU_{0})<\infty$.

Consider two processes with jumps \ce \left\{\begin{array}{ll}
\eta_t=\eta_0+\int^t_0e_s\,\dif s+\int^t_0f_s\,\dif W_s+\int^{t+}_0\int_{\mU_{0}}g(s-,u)\,\tilde{N}_{q}(\dif s,\dif u),\\
\xi_t(x)=\xi_0(x)+\int^t_0E_s(x)\,\dif s+\int^t_0F^l_s(x)\,\dif
W^l_s+\int^{t+}_0\int_{\mU_{0}}G(s-,x,u)\,\tilde{N}_{q}(\dif s,\dif
u),
\end{array}
\right. \de where $\eta_{\cdot},e,f,g(\cdot,u)$ for $u\in\mU_0$ are
predictable processes valued in $\mR^d,\mR^d,\mR^{d\times m},\mR^d$,
respectively, and $\xi_{\cdot}(x),E_{\cdot}(x),F^{l}_{\cdot}(x),
G(\cdot,x,u)$ for $x\in\mR^d,u\in\mU_0$ and $l=1, 2,...,d$ are real
predictable processes; $\{q_t,t\geq 0\}$ is a stationary
$\cF_0^t$-adapted Poisson point process with values in $\mU$ and
characteristic measure $n$. Let $N_{q}((0,t],\dif u)$ be the
counting measure of $q_{t}$ such that $\mE N_{q}((0,t],A)=tn(A)$ for
$A\in\cU$. Denote
$$
\tilde{N}_{q}((0,t],\dif u):=N_{q}((0,t],\dif u)-tn(\dif u),
$$
the compensator of $q_{t}$. Moreover $e,f,g$ and
$E_{\cdot}(x),F^{l}_{\cdot}(x), G(\cdot,x,\cdot)$ satisfy the
following integrable conditions: for $t\in\mR_{+}$ \ce
\int^{t}_0|e_s|\,\dif s<\infty,\int^{t}_0|f_s|^2\,\dif s<\infty,
\int^{t+}_0\int_{\mU_{0}}|g(s-,u)|^2\,n(\dif u)\dif s<\infty,\\
\int^{t}_0|E_s(x)|\,\dif
s<\infty,\sum\limits_l\int^{t}_0|F^l_s(x)|^2\,\dif s<\infty,
\int^{t+}_0\int_{\mU_{0}}|G(s-,x,u)|^2\,n(\dif u)\dif s<\infty. \de

\bp\label{iv} Assume that
$\xi_t(\omega,\cdot)\in\cC_b^2(\mR^{d},\mR)$,
$E_t(\omega,\cdot)\in\cC(\mR^{d},\mR)$,
$F^l_t(\omega,\cdot)\in\cC_b^1(\mR^{d},\mR)$ and
$G(t,\omega,\cdot,u)\in\cC_b^1(\mR^{d},\mR)$ for
$(t,\omega)\in\mR_{+}\times\Omega$ and $u\in\mU_0$. Then \ce
&&\dif\xi_t(\eta_t)\\
&=&\left(\frac{\partial \xi_t(\eta_t)}{\partial
x^l}\right)e^l_t\,\dif t +\left(\frac{\partial
\xi_t(\eta_t)}{\partial x^l}\right)f^{lj}_{t}\,\dif W^j_t
+\frac{1}{2}\left(\frac{\partial^2 \xi_t(\eta_t)}{\partial x^l\partial x^i}\right)f^{lj}_tf^{ij}_{t}\,\dif t\\
&&+\int_{\mU_0}\left[\xi_{t-}(\eta_{t-}+g(t-,u))-\xi_{t-}(\eta_{t-})
-\left(\frac{\partial \xi_{t-}(\eta_{t-})}{\partial x^l}\right)g^l(t-,u)\right]\,n(\dif u)\dif t\\
&&+E_t(\eta_t)\,\dif t+F^l_t(\eta_t)\,\dif W^l_t
+f_t^{li}\left(\frac{\partial F^i_t(\eta_t)}{\partial x^l}\right)\,\dif t\\
&&+\int_{\mU_0}\left[\xi_{t-}(\eta_{t-}+g(t-,u))-\xi_{t-}(\eta_{t-})
+G(t-,\eta_{t-}+g(t-,u),u)\right]\,\tilde{N}_q(\dif t,\dif u)\\
&&+\int_{\mU_0}\left[G(t-,\eta_{t-}+g(t-,u),u)
-G(t-,\eta_{t-},u)\right]\,n(\dif u)\dif t.
\de
\ep

The It\^o-Ventzell formula for $t<0$ could be proved by the similar method to one in \cite{q1}.

\bigskip

{\bf Proof of formula (\ref{form}) in Section \ref{itosde}:}\\
 We examine  the right hand side of the formula (\ref{form}), and show that it is equal to the left hand side. By the Fubini
theorem   \cite{pz}, we obtain that
\ce
&&\int_s^t\left(\int_s^ue^{r-u}\frac{\partial\sigma_i}{\partial
x}(h_r^{x,u})D(r,u)\dif \hat{L}^i_r\right)\dif u
-\int_s^t\left(\int_s^ue^{r-u}\sigma_i(h_r^{x,u})\dif \hat{L}^i_r\right)\dif u\\
&=&\int_s^t\left(\int_s^u\frac{\partial}{\partial v}\left[e^{r-v}\sigma_i(h_r^{x,v})\right]|_{v=u}\dif \hat{L}^i_r\right)\dif u\\
&=&\int_s^t\left(\int_s^t\frac{\partial}{\partial v}\left[e^{r-v}\sigma_i(h_r^{x,v})\right]|_{v=u}I_{[s,u]}(r)\dif \hat{L}^i_r\right)\dif u\\
&=&\int_s^t\left(\int_s^t\frac{\partial}{\partial v}\left[e^{r-v}\sigma_i(h_r^{x,v})\right]|_{v=u}I_{[s,u]}(r)\dif u\right)\dif \hat{L}^i_r\\
&=&\int_s^t\left(\int_s^r\frac{\partial}{\partial v}\left[e^{r-v}\sigma_i(h_r^{x,v})\right]|_{v=u}I_{[s,u]}(r)\dif u\right)\dif \hat{L}^i_r\\
&&+\int_s^t\left(\int_r^t\frac{\partial}{\partial v}\left[e^{r-v}\sigma_i(h_r^{x,v})\right]|_{v=u}I_{[s,u]}(r)\dif u\right)\dif \hat{L}^i_r\\
&=&e^{-t}\int_s^te^r\sigma_i(h_r^{x,t})\dif
\hat{L}^i_r-\int_s^t\sigma_i(h_r^{x,r})\dif \hat{L}^i_r.
\de
This proves the formula (\ref{form}) in Section \ref{itosde}.


\begin{thebibliography}{999}

\bibitem{da} D. Applebaum: {\it L\'evy Processes and Stochastic Calculus}. Second Edition, Cambridge Univ. Press, Cambridge, 2009.

\bibitem{la} L. Arnold: {\it Random dynamical systems,} Springer-Verlag
Berlin, 1998.


\bibitem{Billingsley}  P. Billingsley, {\em Weak Convergence of Probability Measures},  John Wiley/Sons, New York, Second Edition, 1999.

\bibitem{tmr} E. A. Coayla-Teran, S. E. A. Mohammed and P. R. C. Ruffino: Hartman-Grobman
theorem along hyperbolic stationary trajectories, {\it Discrete and
Continuous Dynamical Systems,} 17(2007)281-292.

\bibitem{fl} F. Flandoli and H. Lisei: Stationary conjugation of flows for
parabolic SPDEs with multiplicative noise and some applications,
{\it Stochastic Analysis and Applications,} 22(2005)1385-1420.

\bibitem{fs} F. Flandoli and B. Schmalfuss: Random attractors for the 3D stochastic
Navier Stokes equation with multiplicative noise, {\it Stochastics
and Stochastic Rep.,} 59(1996)21-45.


\bibitem{fk} T. Fujiwara, H. Kunita: Stochastic differential equations of jump type and
L\'evy processes in diffeomorphisms group, {\it J. Math. Kyoto
Univ.,} 25(1985)71-106.

\bibitem{g} M. Grigoriu: Lyapunov exponents for nonlinear systems with Poisson white noise,
{\it Physics Letters A,} 217(1996)258-262.

\bibitem{hwy} S. He, J. Wang and J. Yan: {\it Semimartingale Theory and
Stochastic Calculus,} Science Press and CRC Press Inc., 1992.

\bibitem{iw} N. Ikeda, S. Watanabe: {\it Stochastic differential equations
and diffusion processes,} 2nd ed., North-Holland/Kodanska,
Amsterdam/Tokyo, 1989.


\bibitem{il2} P. Imkeller and C. Lederer: The cohomology of stochastic and random differential
equations, and local linearization of stochastic flows, {\it
Stochastics and Dynamics,} 2(2002)131-159.

\bibitem{is} P. Imkeller and B. Schmalfuss: The conjugacy of stochastic and random differential
equations and the existence of global attractors, {\it Journal of
Dynamics and Differential Equations,} 13(2001)215-249.

\bibitem{ks} G. Kager and M. K. R. Scheutzow: Generation of one-sided
random dynamical systems by stochastic differential equations, {\it
Electronic Journal of Probability,} 2(1997)1-17.

\bibitem{kuni} H. Kunita: Stochastic differential equations based on L\'evy
processes and stochastic flows of diffeomorphisms. In {\it Real and
Stochastic Analysis: New Perspectives}, M. M. Rao (eds.), Birkh\"auser
Boston, 305-373, 2004.

\bibitem{lb} C. W. Li and G. L. Blankenship: Almost sure stability of linear
stochastic systems with Poisson process coefficients, {\it SIAM J.
Appl. Math.,} 46(1986)875-911.

\bibitem{ma} S. I. Marcus: Modelling and approximations of stochastic differential equations driven by
semimartingales, {\it Stochastics,} 4(1981)223-245.

\bibitem{pz} S. Peszat and J. Zabczyk: {\it Stochastic Partial Differential Equations
with L\'evy Noise,} Cambridge University Press, 2007.

\bibitem{q1} H. Qiao: Homeomorphism Flows for Non-Lipschitz SDEs Driven by L\'evy
Processes, {\it Acta Mathematica Scientia,} 32(2012)1115-1125.

\bibitem{q2} H. Qiao and J. Duan: Multiplicative ergodic theorem for discontinuous
random dynamical systems, http://arXiv:1204.5050v2.

\bibitem{sa} K. Sato: {\it L\'evy processes and infinitely divisible distributions,}
Cambridge University Press, 1999.

\bibitem{sh} K. R.Schenk-Hopp\'e: Bifurcation scenarios of the noisy Duffing-van der Pol
oscillator, {\it Nonlin. Dynam.,} 11(1996)255-274.

\bibitem{sm} M. Spivak: {\it A Comprehensive Introduction to Differential Geometry,} Vol. 1, 2nd ed.,
Pub. Berkeley, 1979.

\bibitem{Wanner2} T. Wanner:
Linearization of Random Dynamical Systems. In {\em Dynamics
Report},Volumn 4. Spring-Verlog, New York, 1995.


\bibitem{Woy} W. A. Woyczynski: L\'evy processes in the physical
sciences. In \emph{L\'evy Processes: Theory and Applications}, O.
E. Barndorff-Nielsen, T. Mikosch and S. I. Resnick (Eds.),
241-266, Birkh\"auser, Boston, 2001.

\bibitem{x} Y. Xie:
The random attractors of stochastic Duffing-Van der Pol equations with
jumps, {\it Chinese Journal of Applied Probability and Statistics,} 26(2010)9-23.

\end{thebibliography}
\end{document}